\documentclass{article} 
\newcommand{\be}{\begin{equation}}
\newcommand{\ee}{\end{equation}}
\newcommand{\bd}{\begin{displaymath}}
\newcommand{\ed}{\end{displaymath}}
\author{G.Bergdolt}
\title{An algebraic approach to representations of the permutation group.} 
\begin{document} 
\maketitle 
\section{Abstract}
The group algebra of the permutation group is spanned by a set of elements
called projectors.The coordinates of permutations expanded in projectors
are matrix elements of irreducible representations.The projectors for the 
permutation group are products of a Young symmetriser and an antisymmetriser. 
They form non-orthogonal bases of right and left modules. 
The non-orthogonality is compensated by a constant matrix. It turns out 
that this reduces the matrix 
entries to $\{-1,0,+1\}$ . An algorithm to compute the projectors is given.
 \par 
\section{Introduction}
 The group of permutations of n objects is called the symmetric group and is
denoted $S_n$. A permutation group is a subgroup of $S_n$. Since any 
symmetric group is a subgroup of a symmetric group with a higher number of
symbols we can state: 'A symmetric group is a permutation group'. The 
converse is not true, since a permutation group is the direct product of
symmetric groups relative to disjoint subsets of symbols. In the following
the names permutation group and subgroup of a permutation group are used for
$S_n$ and its subgroups \par
The permutation group is a classical subject, 
the amount of work devoted to the permutation group can be 
infered from the list of references (several hunderts )
in the review article of James and
Kerber in the collection Encyclopedia of Mathematics $^{1)}$.
The origin of the present paper is the chapter entitled 'Symmetric groups' 
  in Littlewood's book $^{2)}$. It was 
  realized that algebra elements introduced 
  there can be promoted to a basic concept of representation theory. 
These elements are called projectors below 
since the relations they satisfy generalise the idempotency
 of projectors. The
projectors for the permutation group, defined as the product of a symmetriser 
and an antisymmetriser, yield a non-orthogonal
 matrix representation. The homomorphy of the matrix representation can also
be obtained if a constant matrix is intersped in the matrix product. 
 It is known at least since Wigner $^{3)}$ that the
representation matrices of a finite group can be taken unitary. 
 Unitarity requires complex numbers for the matrix entries. In the case of 
the permutation group, a property called ambivalence (an element and its inverse
are in the same conjugation class) entails that the representation matrices
can be chosen real orthogonal$^{ 4)}$. Orthogonality still involves square roots 
and hence irrational numbers. 
Dropping orthogonality, the matrix entries can be integers. 
When taking the non-orthogonality
 into account by a constant matrix, the entries are
  reduced to $\{-1,0,+1\}$. This set is the group $Z_3$ under
addition modulo 3. Whether this fact is of significance or not is 
an open question. \par
 The projector formalism is described in Section 3. In Section 4 
 it is shown that the product relations defining projectors are satisfied by
 the products of a Young symmetriser and a Young antisymmetriser. The proof
 that matrix entries can be reduced to $\{-1,0,+1\}$ is given in Section 6.
 Units and characters are examined in Section 7. Conclusions are in Section 8. 
 An algorithm for determining the 
 coordinates of projectors is given in Appendix A. \par
 \section{ Projectors}
 Let A be the algebra spanned by $n^2$ elements $p_{ij}$ called projectors
  satisfying the relations:
 \be p_{ij}p_{kl}=g_{jk}p_{il} \ee
 where $g_{jk}\in k$ is a scalar, k denotes a field.\par
 Let $X,Y \in A$ be defined by coordinates $X=x^{ij}p_{ij}$ and $Y=y^{ij}
 p_{ij}$, 
 where  the summation convention is used.
 The coordinates of $Z=XY$ are given by:
 $$ z^{il}=x^{ij}g_{jk}y^{kl}.$$ 
 The algebra A is seen to be a matrix algebra where the product $\circ$ is
  defined by
 a constant matrix $g$. $$ x \circ y=xgy .$$
 Note that by multiplying to the right or left by the matrix $g$ a conventional
 matrix algebra is obtained. In matrix form $(zg)=(xg)(yg)$.\par
  The direct sum of several such algebras is spanned by projectors $p_{\lambda 
  ij}$ where $\lambda $ labels the subalgebras. 
  The defining relations satisfied by the
 projectors are:
  \be  p_{\lambda ij}p_{\mu kl}= 
 \delta_{\lambda \mu}g_{\lambda jk}p_{\lambda il}.\ee 
 It follows from these relations that \\
 i) the projectors $p_{\lambda ij}$ with $\lambda$ and $i$ fixed span a right
    ideal,\\
 ii) the projectors with $\lambda$ and $j$ fixed span a left ideal,\\
 iii) projectors with $\lambda $ fixed span a subalgebra. \par
 Assume that $\lambda $ ranges over all irreducible representations of a finite
 group: the number of projectors is equal to the number of group elements, i.e.
 the order of the group. This follows from the theorem:
' A regular representation
 contains each irreducible representation with a multiplicity equal to its 
 dimension.'
 If $m_\lambda $ is the dimension of the irreducible representation
 and $n$ the order of the group, we have $$ n = \sum_\lambda 
   m^2_\lambda .$$ The right-hand side is the number of projectors. Note that 
   the representation obtained by the direct sum of irreducible representations
   has dimension $$\sum_\lambda m_\lambda < n,$$ 
 and thus is not the regular representation.  
 \par 
 Since a finite group can be reconstructed from the set of its irreducible
 representations, 
  the equality of the number of projectors and group elements entails 
   that the projectors form a basis of the group algebra.\par   
  Elements of the group algebra $kS_n$ satisfying these relations are
 defined in the next section. \par
 \section{Projectors for the permutation group}
 The representation theory of the permutation group can be sketched as follows:
 Irreducible representations of $S_n$ are characterised by partitions of $n$. A
 partition defines a Young frame. Young frames filled with $n$ symbols are Young
 tableaus. A lexical order of the symbols is adopted. A Young frame with the
 symbols in rows and columns in lexical order is a standart Young tableau.
 The number of standard Young tableaus for a 
 given partition is the dimension of the irreducible
 representation. 
  A Young tableau is defined by a partition and the sequence of symbols obtained
 by reading the tableau from left to right, top to bottom. The sequences are
 ordered by the first differing symbol in the sequences. 
 This yields an ordering of Young tableaus.\par
Let $i$ label Young tableaus and let $P_i$ be the row symmetriser 
 of the tableau $T_i$. i.e. the
 sum of all permutations permuting symbols in rows of $T_i$.
  Let $N_i$  be the
 column antisymmetriser i.e. the sum of signed permutations permuting the 
 symbols in columns of $T_i$ with sign +  for even, - for odd permutations.
 Let $\sigma_{ij}$ be the permutation which permutes 
 the sequence of tableau $T_i$
 into the sequence of tableau $T_j$. We have $N_i\sigma_{ij}=\sigma_{ij}N_j$ 
 and $P_i \sigma_{ij} = \sigma_{ij} P_j.$ i.e. $\sigma_{ij}$ is an intertwiner.
 \\
 Define $p_{ij} \in kS_n$ by:
 \be p_{ij}=P_i \sigma_{ij}N_j. \ee 
 Proposition: The elements $p_{ij}$ satify the projector relations (1).\par
 Proof: 
 The set of permutations permuting the symbols in columns of a Young tableau
 form a subgroup of $S_n$ as do the permutations permuting the symbols in 
 rows.
  Denote these subgroups also by $N_j$ and $P_j$ for tableau $T_j$.
   For a 
 permutation $v_j \in N_j$ we have $N_j v_j = \pi (v_j )N_j$
 where $\pi (v_j )$ is the parity of $v_j $. For a permutation
 $h_k \in P_k$ we have $ h_k P_k = P_k$.\\
 If an odd permutation is contained in $N_j$ and in $P_k$ we have 
 $N_j P_k = - N_j P_k $ and the product is null. This is the case if two 
 symbols are in the same column in $T_j$ and the 
 same row in $T_k$. The transposition
 of the two symbols is an odd permutation contained in $N_j$ and $P_k$.
 If no permutation except the identity is contained in $N_j$ and $P_k$ a
 permutation $ v_j \in N_j$ can be found wich puts a symbol of $T_j$
 in the same row as in $T_k$. A permutation $ h_k \in P_k$ then exists
 which puts the symbols in the same position as in $T_k$, hence $\sigma_{jk}=
  v_j h_k $.\\It follows that:
 \begin{eqnarray*}
  p_{ij}p_{kl} &=& P_i \sigma_{ij} N_j P_k \sigma_{kl} N_l \\
               &=&  \pi(v_j)(P_iN_i)^2 \sigma_{il}. \end{eqnarray*}
     From Littlewood $^{2)}$ we have the result $(P_i N_i)^2= ( n!/m) P_iN_i$ 
     where $m$ is the number of standart Young tableaus. Relation (1) follows 
 with \bd  g_{jk}= \frac {n!}{m} \pi(v_j).\ed \par 
 If the Young tableaus are ordered according to the ordering defined above then
 $N_i P_j =0$ if $T_i < T_j$ and the $g$ matrix is lower triangular.
  The proof of the proposition ' If $T_i<T_j$ the first differing symbols
 are in the same column in $T_i$ and the
 same row in $T_j$' is given by Littlewood $^{2)}$. As shown above this
 entails that $N_iP_j=0$. \par  
 Let $\lambda,\mu $ label partitions and let $N^\lambda_i$, $P^\mu_i$
 denote the corresponding antisymmetriser and symmetriser.
  According to another result of
 Littlewood$^{ 2)}$ $(P^\lambda_i N^\lambda_i)(P^\mu_j N^\mu_j)=0 $
 if $\lambda \ne
 \mu$. Hence the general relation (2) is satisfied by the projectors. \\
 Remark: Another set of projectors is provided by $p^*_{ij}=N_i \sigma_{ij} 
 P_j$ . The sets $\{ p_{ij}\}$ and $\{p^*_{ij}\}$ are related by the 
 involution $s \rightarrow s^{-1}$ of $S_n$.
 \section{Matrix representations of $S_n$} 
  Denote by $s_a, a=1, \dots,n!$ the permutations of $S_n$. Expanding $s_a$ in
  the projector basis we have   
  \be s_a=x_a^{\lambda ij}p_{\lambda ij}.\ee
  where $x_a^{\lambda ij}$ are coordinates and the summation is over $\lambda i
  j$.
   From the product relations (2) it
  follows that the coordinates with $\lambda$ fixed provide an irreducible 
  matrix representation of $S_n$. If $s_c=s_as_b$ then 
  $$x_c^{\lambda il}=x_a^{\lambda ij} g_{\lambda jk} x_b^{\lambda kl} $$
  where the summation is over $j$ and $k$. \par 
   The relation (3)
   defines the projectors as linear combinations 
  of permutations  Hence: \be p_{\lambda ij} = y_{\lambda ij}^a s_a \ee 
   where the summation is over $a$. \\
  Note that the coordinates $y^a_{\lambda ij}$ have values $\{-1,0,+1\}$.\\
  Proof: The subgroups $N_i$ and $P_j$ have no permutation in common exept the
  Identity $e$. If $c',c'' \in N_i$ and $r',r'' \in P_i$ are such that  
$c'r'= c''r''$ then
  since $r''^{-1}r' \in N_i$ and $c''c'^{-1}\in P_i$ we have $r'=r''$ and 
$c'=c''$.\par
   The linear system (4) is the inverse of (5) ; it follows that the 
matrices of 
   coordinates are inverses $(x)=(y)^{-1}$.
   The matrix of coordinates $(y)$ is given by the definition of projectors.
  The problem of determining the irreducile representations 
of $S_n$ is then solved in principle. \\
  Example: The group $S_3$\\ Define the permutations by the sequence obtained 
  applying the permutation to $[123]$ : \par
  $ s_1=[123], s_2=[132], s_3=[213],
  s_4=[231], s_5=[312], s_6=[321] $\\ 
  The partitions are labeled :1=(3), 2=(2,1), 3=(111).\\ The definition of 
  projectors yields the linear system:
  \bd \left( \begin{array}{c} p_{111} \\ p_{211} \\ p_{212} \\ p_{221} \\
  p_{222} \\ p_{311} \end{array} \right) = \left( \begin{array}{rrrrrr}
   1 &  1 &  1 &  1 &  1 &  1 \\
   1 &  0 &  1 &  0 & -1 & -1 \\
   0 &  1 & -1 & -1 &  1 &  0 \\
   0 &  1 &  0 &  1 & -1 & -1 \\
   1 &  0 & -1 & -1 &  0 &  1 \\
   1 & -1 & -1 &  1 &  1 & -1 \end{array} \right) \left( \begin{array}{c} s_1 \\
  s_2 \\ s_3 \\ s_4 \\ s_5 \\ s_6 \end{array} \right) \ed 
  The projector relations (2) are satisfied with a $g$ matrix $g=
    diag (6,3,3,3,3,6)$. The linear system can be inverted ; permutations are 
    related to projectors by: 
   \bd \left( \begin{array}{c} s_1 \\ s_2 \\ s_3 \\ s_4 \\ s_5 \\ s_6 
   \end{array} \right) = \left( \begin{array}{rrrrrr}
   1 &  1 &  0 &  0 &  1 &  1 \\
   1 &  0 &  1 &  1 &  0 & -1 \\
   1 & -1 &  1 & -1 &  0 &  1 \\
   1 &  1 & -1 &  0 & -1 & -1 \\
   1 &  0 & -1 &  1 & -1 &  1 \\ 
   1 & -1 &  0 & -1 &  1 & -1 \end{array} \right) \left( \begin{array}{c}
    p_{111}/6 \\ p_{211}/3 \\ p_{212}/3 \\ p_{221}/3 \\ p_{222}/3 \\
   p_{131}/6 \end{array} \right) \ed 
   As pointed out the entries of the representation matrices are the 
   coordinates of the permutations. Moreover the normalisation factors of
   the projectors and the entries of the $g$ matrix cancel so that the two
   dimensional representation can be read off the matrix above:
   \bd 
   s_1 = \left( \begin{array}{rr} 1 & 0 \\ 0 & 1 \end{array} \right) ,
   s_2 = \left( \begin{array}{rr} 0 & 1 \\ 1 & 0 \end{array} \right) ,
   s_3 = \left( \begin{array}{rr} 1 & 0 \\-1 &-1 \end{array} \right), 
   \ed
  \bd 
   s_4 = \left( \begin{array}{rr}-1 &-1 \\ 1 & 0 \end{array} \right) ,
   s_5 = \left( \begin{array}{rr} 0 & 1 \\-1 &-1 \end{array} \right) ,
   s_6 = \left( \begin{array}{rr}-1 &-1 \\ 0 & 1 \end{array} \right) .
   \ed 
\section{Reduced entries}
 An algorithm to compute the coordinates of projectors expanded in 
 permutations is described in Appendix A.\\ 
 We have seen that the matrix of coordinates of permutations expanded in
 projectors can be obtained by inversion of the matrix of projector 
 coordinates. There is however a shortcut which avoids the 
 inversion of an $n!$ dimensional matrix. 
 The coordinates $(x)$ and $(y)$ are related by:
 \be y_{\lambda ij}^{b^{-1}}= \frac{m_\lambda }{n!}\sum_{rs}g_{\lambda jr}
 x_b^{\lambda rs}g_{\lambda sr}.\ee
 Here $b$ and $b^{-1}$ are a permutation and its inverse. A proof of relation
 (6) is given in Appendix B. We now prove the claim made in the Introduction:
 Representation matrices exist with entries restricted to $\{-1,0,+1\}.$
 In relation (6) $\lambda $ is fixed. Set $f = \frac {n!}{m_{\lambda }}$.  
 Writing the quantities with indices 
 i,j,r,s as matrices (6) reads : \bd \tilde y(b^{-1}) = f^{-1}g x(b) g \ed
 The integer factor $f$ was incorporated in the $g$ matrix in Section (4) 
 hence $g'=f^{-1}g$ is a matrix with entries $\{-1,0,+1\}.$
  If projectors are 
 renormalised as $p'_{ij}=f^{-1}p_{ij}$ the relations $p'_{ij}p'_{kl}=g'_{jk}
 p'_{il}$ are satisfied and the coordinates $x'^{ij}=fx^{ij}$ are the matrix
 elements of a representation. Relation (6) is now:
 \bd \tilde y(b^{-1})=g'x'(b)g' \ed 
 If $x'(b)$ defines a representation with matrix $g'$ then $g'x'(b)g'$
 defines a representation with matrix $g'^{-1}$.\\
 Proof: $(g'x'(b)g')g'^{-1}(g'x'(c)g')=g'(x(b)g'x'(c))g'=(g'x'(bc)g').$ 
 Note that the transposed matrix $\tilde y(b^{-1})$ on the left-hand side of
 (6) is a matrix with entries $\{-1,0,+1\}$ 
 so that the claim is proved. 
 The matrix $g'^{-1}$ is the inverse of a lower triangular matrix with reduced
 entries. It follows that the entries of the matrix $g'^{-1}_n$
are integers but not that they are 
 reduced to $\{-1,0,+1\}$.
 \section{ Units and conjugation classes} 
 Let $C_\rho $ be the sum of the permutations in a conjugation class. Recall
 that the number of conjugation classes of a finite dimensional group is 
 equal to the number of irreducible representations. $C_\rho $ commutes with
 all elements of the group algebra. \\
 Define: $$U_\lambda = \sum_{ij} g^{-1}_{\lambda ij} p_{\lambda ij}. $$
 From the product relations (3) it follows that  $U_\lambda \ p_{\mu kl}=
 \delta_{\lambda \mu} p_{\lambda kl}$ and \\$p_{\mu kl} \ U_{\lambda}=
 \delta_{\mu \lambda}p_{\lambda kl}.$ Further $(U_\lambda)^2=U_\lambda $.
 $U_\lambda$ is idempotent,  leaves elements of the subalgebra $\lambda $
 invariant and annihilates all others, i.e. $U_\lambda $ is a unit of the 
 subalgebra $\lambda $. $U_\lambda$ commutes with all elements of the algebra,
 hence $U_\lambda$ is a sum of entire conjugation classes:
 \be U_\lambda = \sum_\rho \chi^\rho_\lambda C_\rho.\ee
 The coefficients are group characters of $S_n$. A permutation cannot be in
 several classes and the elements of $g^{-1}_\lambda $ are integers. It follows
 that the characters defined above are integers. \\
 Example: The group $S_3$.\par 
 Permutations in cycle notation are denoted by$(\ldots )$ and in sequence
 notation by $[\ldots]$.
 The conjugation classes are:
\begin{eqnarray*} C_1 &=& (1)(2)(3) ,\\
                 C_2 &=& (12)+(23)+(13) ,\\
		 C_3 &=& (123)+(132) . \end{eqnarray*}
  The units are given by: 
 \begin{eqnarray*} U_1 &=& [123]+[132]+[312]+[213]+[231]+[321],\\
     U_2 &=& 2[123]-[312]-[231],\\
     U_3 &=& [123]-[132]+[312]-[213]+[231]-[321]. \end{eqnarray*} 
  The linear system (7) is then:
  \bd \left( \begin{array}{r} U_1 \\ U_2 \\ U_3 \end{array} \right) = 
  \left( \begin{array}{rrr} 1 & 1 & 1 \\ 2 & 0 & -1 \\ 1 & -1 & 1 
  \end{array} \right) \left( \begin{array}{r} C_1 \\ C_2 \\ C_3 \end{array}
  \right) . \ed 
  \section{Conclusion}
   The approach is resticted to finite groups. 
  With respect to the usual representation theory the approach 
  described above could be called natural. Orthonormalisation requires 
  irrational numbers as matrix entries. In the $g$ matrix scheme entries
  are integers restricted to $\{-1,0,+1\}$.
.For example the set of matrices with reduced entries given in 
  section 5 is closed under matrix multiplication. The
  $g$ matrix defines a vector space isomorphism between right and left
modules and denotes an intrinsic property of the $kS_n$ algebra. The algorithm
described in Appendix B has been implemented in a computer program.For $n<4$
the $g$ matrices are diagonal. The first non-diagonal $g$
matrix occurs for $n=5$ and partition (3,2). The first matrix $g'^{-1}$
with entries not in $\{-1,0,+1\}$  occurs for$n=7$ and partition $(3,2,2)$.
The representations matrices are always matrices with reduced entries.\par
The reduced entries feature has not found applications to physics up to now.
 The theory of
  coherent states in quantum mechanics shows that orthonormalisation is not 
  an unavoidable feature. \par
  \par 
 
\par \par
\appendix{Appendix A}      		  
According to definition (3) $p_{ij}$ is a sum of terms of the form $\pi(v_i) 
  h_j \sigma_{ji}v_i$ 
 where $h_j$ is a row permutation of $T_j$, $v_i$ a column permutation of 
$T_i$ and $ \pi(v_i) = +1$ for even, $-1$ for odd permutations.
 If a permutation $s$ can be factorised as $s=h_j\sigma_{ji}v_i$ the coordinate 
 of $p_{ji}$ along $s$ is given by $y^s_{ji}=\pi (v_i)$. \\
 Define the action of a permutation on a Young tableau as a right action:
 $T_i \rightarrow T_i s = T_s$ The conjugates of the permutations $v_i$ by
$ s $ : $s v_i s^{-1}=v_s$ form a subgroup ; 
 the subgroup of column permutations of $T_s$.
 Note that the parity of $v_i$ and $v_s$ are the same. We have 
 $$ v_i^{-1} s = \sigma_{ij} h_j = s v_s^{-1}. $$  
 If this permutation
 acts on $T_i$ we obtain $T_s v_s^{-1} = T_j h_j $. A column permutation
 of $T_s$ exists which puts the symbols in the same rows as in $T_j$. This 
 yields the following algorithm to compute the coordinates of a projector. \par
 Determine the non-standart Young tableau $T_s=T_i s$. 
 Set to 0 a parameter $k$.
 Scan the columns of the non standard Young tableau $T_s$ 
 and take the sequence of symbols in a column. 
 Locate the same symbols in $T_j$ ordered according to the
 rows they occupy. This is not possible if two such symbols are in the same row 
 and in that case the algorithm ends with a zero coordinate. 
 Otherwise the two sequences of symbols define a permutation.
 Determine the cycle structure of the permutation.  For a cycle of
 length $l$ add $l-1$ to the parameter $k$. The coordinate is given by:
 $y_{ij}^s = 1-2(k \  mod \ 2)$\\
 Note that the matrix elements of $g$ are obtained by setting $s=e$. 
 \appendix{Appendix B} 
 Since the (x) and (y) matrices are inverses we have 
 $$ \sum_b y^b_{\lambda ij} x_b^{\mu kl} =\delta^\mu_\lambda \delta^k_i
 \delta_j^l ,$$ 
 $$ \sum_{\lambda ij} x^{\lambda ij}_a y^b_{\lambda ij} = \delta^a_b.$$ 
 where $a,b$ are permutations. Recall the representation relation 
 $$ x_c^{\lambda il} = \sum_{jk} x_a^{\lambda ij}g_{\lambda jk}x_b^{\lambda
    kl},$$ where $c=ab$. Now \\
 i) multiply by $y^a_{\lambda ij}$ and sum over a,\\
 ii)multiply by $y^c_{\mu ks}$ and sum over $\mu ks$,\\
 the result is :
 \be y_{\lambda ij}^{cb^{-1}}=\sum_{rs}x_b^{\lambda rs}g_{\lambda jr} 
    y^c_{\lambda is}. \ee
 If $c=b$ the relation is:
 $$y^e_{\lambda ij}=\sum_{rs}x_b^{\lambda rs}g_{\lambda jr}y^b_{\lambda is}.
 $$  Summing over $b$ gives :
 $$ n! y^e_{\lambda ij}=m_\lambda g_{\lambda ji}.$$
 If in (8) $c=e$ we have :
 $$y^{b^{-1}}_{\lambda ij}=\sum_{rs} g_{\lambda jr} x_b^{\lambda rs} 
 y^e_{\lambda is}.$$ From the last two relations, relation (6) follows.
 
\end{document}